\documentclass[article]{biometrika}

\usepackage{times}
\usepackage{bm}
\usepackage{natbib}
\usepackage{epstopdf} 

\usepackage{amsmath}
\usepackage{bm,bbm}
\usepackage{amsfonts}
\usepackage{graphics}
\usepackage{graphicx}
\usepackage{caption,subcaption}

\newcommand{\R}{\mathbb{R}}

\newcommand{\bA}{A}
\newcommand{\bB}{B}

\newcommand{\bE}{E}

\newcommand{\bH}{H}
\newcommand{\bI}{I}

\newcommand{\bM}{M}

\newcommand{\bP}{P}

\newcommand{\bR}{R}

\newcommand{\bT}{T}
\newcommand{\bU}{U}

\newcommand{\bW}{W}
\newcommand{\bX}{X}
\newcommand{\bY}{Y}

\newcommand{\bSig}{\Sigma}
\newcommand{\bLam}{\Lambda}
\newcommand{\bXi}{\Xi}
\newcommand{\bGamma}{\Gamma}

\newcommand{\rmF}{\textrm{F}}

\newcommand{\rmRank}{\textrm{rank}}

\newcommand{\rmMax}{\textrm{max}}
\newcommand{\rmMin}{\textrm{min}}

\newcommand{\bZero}{0}

\newcommand{\Prob}{\mathbb{P}}

\newcommand{\ttinf}{2\rightarrow\infty}

\newcommand{\defO}{\mathbb{O}}

\newcommand{\probO}{O_{\mathbb{P}}}

\newtheorem{question}[theorem]{Question}

\begin{document}

\jname{Submitted to Biometrika}
\jyear{NA}
\jvol{NA}
\jnum{NA}


\markboth{J.~Cape \and M.~Tang \and C.~E.~Priebe}{Eigenvector deviations and fluctuations}

\title{Signal-plus-noise matrix models:\\ eigenvector deviations and fluctuations}

\author{J. CAPE, M. TANG, \and C. E. PRIEBE}
\affil{
	Department of Applied Mathematics and Statistics, Johns Hopkins University\\
	3400 North Charles Street; Baltimore, Maryland 21218 U.S.A. \email{joshua.cape@jhu.edu}
	\email{minh@jhu.edu}
	\email{cep@jhu.edu}}

\maketitle

\begin{abstract}
Estimating eigenvectors and low-dimensional subspaces is of central importance for numerous problems in statistics, computer science, and applied mathematics. This paper characterizes the behavior of perturbed eigenvectors for a range of signal-plus-noise matrix models encountered in both statistical and random matrix theoretic settings. We prove both first-order approximation results (i.e.~sharp deviations) as well as second-order distributional limit theory (i.e.~fluctuations). The concise methodology considered in this paper synthesizes tools rooted in two core concepts, namely (i)~deterministic decompositions of matrix perturbations and (ii)~probabilistic matrix concentration phenomena. We illustrate our theoretical results via simulation examples involving stochastic block model random graphs.
\end{abstract}


\section{Introduction}
\label{sec:Intro}
This paper considers the setting where $\bM$ and $\bE$ are large $n \times n$ symmetric real-valued matrices with $\hat{\bM}=\bM+\bE$ representing an additive perturbation of $\bM$ by $\bE$. For $n \times r$ matrices $\bU$ and $\hat{\bU}$ whose columns are orthonormal eigenvectors corresponding to the $r \ll n$ leading eigenvalues of $\bM$ and $\hat{\bM}$, respectively, we ask:
\begin{question}[1]
	\label{q:main}
	How entrywise close are the matrices of eigenvectors $\bU$ and $\hat{\bU}$?
\end{question}
Under quite general structural assumptions on $\bU$, $\bM$, and $\bE$, our main results address Question~\ref{q:main} both at the level of first-order deviations and at the level of second-order fluctuations. Theorems~\ref{thrm:firstOrder}~and~\ref{thrm:firstOrderExtension} quantify the entrywise closeness of $\hat{\bU}$ to $\bU$ modulo a necessary orthogonal transformation $\bW$ which will subsequently be made precise. Theorem~\ref{thrm:secondOrder} states a multivariate distributional limit result for the rows of the matrix $\hat{\bU}-\bU\bW$.

Numerous problems in statistics consider the eigenstructure of large symmetric matrices. Prominent examples of such problems include (spike) population and covariance matrix estimation \citep{silverstein1984some,silverstein1989eigenvectors,johnstone2001PCA,yu2014useful} as well as principal component analysis \citep{jolliffe1986principal,nadler2008finite,paul2007asymptotics}, problems which have received additional attention and windfall as a result of advances in random matrix theory \citep{bai2010spectral,benaych2011eigenvalues,paul2014random}. Within the study of networks, the problem of community detection and success of spectral clustering methodologies have also led to widespread interest in understanding spectral perturbations of large matrices, in particular graph Laplacian and adjacency matrices \citep{rohe2011spectral,lei2015consistency,sarkar2015role,le2017concentration,tang2018limit}. Towards these ends, recent ongoing and concurrent efforts in the statistics, computer science, and mathematics communities have been devoted to obtaining precise entrywise bounds on eigenvector perturbations \citep{fan2018ell, cape2018two, eldridge2018unperturbed}. See also \cite{mao2017estimating}, \cite{abbe2017entrywise}, and \cite{tang2017asymptotically}.

This paper distinguishes itself from the literature by presenting both deviation and fluctuation results within a concise yet flexible signal-plus-noise matrix model framework amenable to statistical applications. We extend the machinery and perturbation considerations introduced in \cite{cape2018two} in order to obtain strong first-order bounds. We then demonstrate how careful analysis within a unified framework leads to second-order multivariate distributional limit theory. Our characterization of eigenvector perturbations relies upon a matrix perturbation series expansion together with an approximate commutativity argument for certain matrix products.

The results in this paper apply to principal component analysis in spike matrix models, including those of the form $\bY = \lambda u u^{\top} + n^{-1/2}\bE$ where $u\in\R^{n}$ denotes a spike (signal) unit vector and $\bE\in\R^{n \times n}$ denotes a random symmetric (noise) matrix. We consider the super-critical regime, $\lambda~>~1$, for which it is known, for example, that the leading eigenvector $\hat{u}$ of $\bY$ has non-trivial correlation with $u$ when $\bE$ is drawn from the Gaussian orthogonal ensemble, namely $|\langle \hat{u}, u \rangle|^{2} \rightarrow 1 - 1/\lambda^{2}$ almost surely \citep{benaych2011eigenvalues}. This paper obtains stronger local results for spike estimation in the presence of sufficient eigenvector delocalization and provided the signal in $\lambda~\gg~1$ is sufficiently informative with respect to $\bE$. Loosely speaking, we establish that $\|\hat{u}-u\|_{\infty} \le C(\log n)^{c}\lambda^{-1}\|u\|_{\infty}$ with high probability for some positive constants $C$ and $c$, and we prove that $n(\hat{u}_{i}-u_{i})$ is asymptotically normally distributed. Our results hold more generally for $r$-dimensional spike models exhibiting eigenvalue multiplicity and for $\bE$ exhibiting a heterogeneous variance profile.

\section{Preliminaries}

For $n \times r$ real matrices with orthonormal columns, denoted by $\hat{\bU}, \bU \in \defO_{n,r}$, the columns of $\hat{\bU}$ and $\bU$ each form orthonormal bases for $r$-dimensional subspaces of $\R^{n}$. The distance between subspaces is commonly defined via the notion of \emph{canonical angles} and the C(osine)-S(ine) matrix decomposition which crucially involve the singular values of the matrix $\bU^{\top}\hat{\bU}$. Specifically, by writing the singular values of $\bU^{\top}\hat{\bU}$ as $\sigma_{1} \ge \sigma_{2} \ge \dots \ge \sigma_{r}$, the $r \times r$ diagonal matrix of canonical angles is defined on the main diagonal as $\Theta(\hat{\bU},\bU)_{ii} = \arccos(\sigma_{i})$ for $i \in [r] = \{1,2,\dots,r\}$ \citep[Section~7.1]{bhatia1997matrix}.

One frequently encounters the entrywise-defined matrix $\sin\Theta(\hat{\bU},\bU)\in\R^{r \times r}$ since for the commonly considered spectral and Frobenius matrix norms, small values of $\|\sin\Theta(\hat{\bU},\bU)\|_{\eta}$ indicate small angular separation (distance) between the subspaces corresponding to $\hat{\bU}$ and $\bU$. Importantly, the canonical angle notion of distance between $r$-dimensional subspaces takes into account basis alignment in the form of right-multiplication by an $r \times r$ orthogonal matrix $\bW \in \defO_{r,r} \equiv \defO_{r}$, and each choice of norm $\eta \in \{\cdot,\rm{F}\}$ satisfies \citep[Lemma 1]{cai2018rate}
\begin{equation*}
	\label{eq:minDistanceOrthogonal}
	\|\sin\Theta(\hat{\bU},\bU)\|_{\eta} \le 
	\underset{\bW \in \defO_{r}}{\textrm{inf}}\|\hat{\bU}-\bU\bW\|_{\eta} \le
	\sqrt{2}\|\sin\Theta(\hat{\bU},\bU)\|_{\eta}.
\end{equation*}
In this paper, we focus on matrices of the form
\begin{equation}
	\hat{\bU}-\bU\bW \in \R^{n \times r},
\end{equation}
but we instead consider the two-to-infinity matrix norm which is defined via the $\ell_{2}$ and $\ell_{\infty}$ vector norms for any matrix $\bT$ as $\|\bT\|_{\ttinf}~=~\textrm{sup}_{\|x\|=1}\|\bT x\|_{\infty}$. The quantity $\|\bT\|_{\ttinf}$ has the convenient interpretation of being the maximum Euclidean norm of the rows of $\bT$ and therefore affords the advantage of being invariant with respect to right-multiplication by orthogonal matrices. Our subsequent analysis will be shown to be particularly meaningful when $\bU$ exhibits low/bounded \emph{coherence} \citep{candes2009exact}, i.e.~when $\bU$ is suitably \emph{delocalized} \citep{rudelson2015delocalization}, in the sense that $\|\bU\|_{\ttinf}$ decays sufficiently quickly in $n$.

For tall, thin matrices $\bT \in \R^{n \times r}$ with $n \gg r$, such as $\hat{\bU}-\bU\bW$, standard norm relations reveal that $\|\bT\|_{\rmMax}~=~\rmMax_{i,j}|\bT_{ij}|$ and $\|\bT\|_{\ttinf}$ differ by at most a factor depending on $r$. The same well-known relationship holds for the spectral and Frobenius norms, $\|\bT\|$ and $\|\bT\|_{\rmF}$, since necessarily $\rmRank(\bT) \le r$. In contrast, $\|\bT\|_{\ttinf}$ may in certain cases be much smaller than $\|\bT\|$ by a factor depending on $n$, summarized as
\begin{equation*}
	\|\bT\|_{\rmMax}
	\overset{r}{\asymp} \|\bT\|_{\ttinf}
	\overset{n}{\ll} \|\bT\|
	\overset{r}{\asymp} \|\bT\|_{\rmF}.
\end{equation*}

Taken together, these properties suggest the appropriateness of the two-to-infinity norm when viewing the rows of $\bT$ as a point cloud of residuals in low-dimensional Euclidean space. We refer the reader to \cite{cape2018two} for a more general discussion of the two-to-infinity norm and for more on statistical applications, including community detection and principal subspace estimation which are of interest here. In the current paper, additional model assumptions and more refined technical analysis yield much stronger results for these applications.

\section{Main Results}
\label{sec:MainResults}

\subsection{Setting}
\label{sec:Setting}

Let $\bM \equiv \bM_{n} \in \R^{n \times n}$ be a symmetric matrix with block spectral decomposition given by
\begin{equation}
	\label{eq:block_spectral_decomp}
	\bM
	\equiv [\bU|\bU_{\perp}][\bLam \oplus \bLam_{\perp}][\bU|\bU_{\perp}]^{\top}
	= \bU\bLam\bU^{\top} + \bU_{\perp}\bLam_{\perp}\bU_{\perp}^{\top},
\end{equation}
where the diagonal matrix $\bLam \in \R^{r \times r}$ contains the $r$ largest-in-magnitude nonzero eigenvalues of $\bM$ with $|\bLam_{11}| \ge \dots \ge |\bLam_{rr}| > 0$, and where $\bU \in \defO_{n,r}$ is an $n \times r$ matrix whose orthonormal columns are the corresponding eigenvectors of $\bM$. The diagonal matrix $\bLam_{\perp} \in \R^{(n-r)\times(n-r)}$ contains the remaining $n-r$ eigenvalues of $\bM$ with the associated matrix of orthonormal eigenvectors $\bU_{\perp} \in \defO_{n,(n-r)}$. Let $\bE \in \R^{n \times n}$ be a symmetric matrix, and write the perturbation of $\bM$ by $\bE$ as
$ \hat{\bM} \equiv \bM + \bE = \hat{\bU}\hat{\bLam}\hat{\bU}^{\top} + \hat{\bU}_{\perp}\hat{\bLam}_{\perp}\hat{\bU}_{\perp}^{\top}$.

\begin{assumption}
	\label{assump:sparsity}
	Let $\rho_{n}$ denote a possibly $n$-dependent scaling parameter such that $(0,1] \ni \rho_{n} \rightarrow c_{\rho} \in [0,1]$ as $n \rightarrow\infty$, with $n\rho_{n} \ge c_{1}(\log n)^{c_{2}}$ for some constants $c_{1}, c_{2} \ge 1$.
\end{assumption}
\begin{assumption}
	\label{assump:evalSignal}
	There exist constants $C,c>0$ such that for all $n \ge n_{0}(C,c)$, $|\bLam_{rr}| \ge c (n \rho_{n})$ and $|\bLam_{11}||\bLam_{rr}|^{-1} \le C$, while $\bLam_{\perp} \equiv \bZero$.
\end{assumption}
\begin{assumption}
	\label{assump:specConcentration}
	There exist constants $C,c >0$ such that $\|\bE\| \le C(n\rho_{n})^{1/2}$ with probability at least $1-n^{-c}$ for all $n \ge n_{0}(C,c)$, written succinctly as $\|\bE\| = O_{\Prob}\{(n\rho_{n})^{1/2}\}$.
\end{assumption}

Assumption~\ref{assump:sparsity} introduces a sparsity scaling factor $\rho_{n}$ for added flexibility. This paper considers the large-$n$ regime and often suppresses the dependence of (sequences of) matrices on $n$ for notational convenience.

Assumption~\ref{assump:evalSignal} specifies the magnitude of the leading eigenvalues corresponding to the leading eigenvectors of interest. For simplicity and specificity, all leading eigenvalues are taken to be of the same prescribed order, and the remaining eigenvalues are assumed to vanish. Remark~\ref{rem:moreAssumption2} briefly addresses the situation when the leading eigenvalues differ in order of magnitude, when $\bLam_{\perp} \neq \bZero$, and when the (spike) dimension $r$ is unknown.

Assumption~\ref{assump:specConcentration} specifies that the random matrix $\bE$ is concentrated in spectral norm in the classical probabilistic sense. Such concentration holds widely for random matrix models where $\bE$ is centered, in which case $\hat{\bM}$ has low rank expectation equal to $\bM$. The advantage of Assumption~\ref{assump:specConcentration} when coupled with Assumption~\ref{assump:evalSignal} is that, together with an application of Weyl's inequality \citep[Corollary~3.2.6]{bhatia1997matrix}, the implicit signal-to-noise ratio terms behave as $\|\bE\||\bLam_{rr}|^{-1}, \|\bE\||\hat{\bLam}_{rr}|^{-1} = \probO\{(n\rho_{n})^{-1/2}\}$. It is straightforward to adapt our analysis and results under less explicit assumptions, albeit at the expense of succinctness and clarity.

Below, Assumption~\ref{assump:entryConcentration} specifies an additional probabilistic concentration requirement that arises in conjunction with the model flexibility introduced via the sparsity scaling factor $\rho_{n}$ in Assumption~\ref{assump:sparsity}. The notation $\lceil \cdot \rceil$ is used to denote the ceiling function.

\begin{assumption}
	\label{assump:entryConcentration}
	There exist constants $C_{\bE}, \nu >0$, $\xi > 1$ such that for all $1 \le k \le k(n) = \lceil\log n / \log (n\rho_{n})\rceil$, for each standard basis vector $e_{i}$, and for each column vector $u$ of $\bU$,
	\begin{equation}
		\label{eq:entryConcentration}
		|\langle\bE^{k}u,e_{i}\rangle|
		\le (C_{\bE} n\rho_{n})^{k/2}(\log n)^{k\xi}\|u\|_{\infty}
	\end{equation}
	with probability at least $1-\exp\{-\nu(\log n)^{\xi}\}$ provided $n \ge n_{0}(C_{\bE}, \nu, \xi)$.
\end{assumption}

Assumption~\ref{assump:entryConcentration} states a higher-order concentration estimate that reflects behavior exhibited by a broad class of random symmetric matrices including Wigner matrices whose entries exhibit subexponential decay and nonidentical variances \citep[modification of Lemma~7.10; Remark~2.4]{erdHos2013spectral}; see also \cite{mao2017estimating}. For example, using our notation, the proof of Lemma~7.10 in \cite{erdHos2013spectral} establishes that $|\langle (C_{\bE} n\rho_{n})^{-k/2}\bE^{k}e,e_{i}\rangle| \le (\log n)^{k \xi}$ with high probability, where $e$ is the vector of all ones and the symmetric matrix $\bE$ has independent mean zero entries with bounded variances. Taking a union bound collectively over $1 \le k \le k(n)$, the standard basis vectors in $\R^{n}$, and the columns of $\bU$, yields an event that holds with probability at least $1-n^{-c}$ for some constant $c > 0$ for sufficiently large $n$.

The function $k(n)$ is fundamentally model-dependent through its connection with the sparsity factor $\rho_{n}$ and satisfies $(n\rho_{n})^{-k(n)/2} \le n^{-1/2}$ for $n$ sufficiently large. In the case when $\rho_{n} \equiv 1$, then $k(n)\equiv 1$, and the behavior reflected in Eq.~(\ref{eq:entryConcentration}) reduces to commonly-encountered Bernstein-type probabilistic concentration. In contrast, when $\rho_{n} \rightarrow 0$ and, for example, $(n\rho_{n}) = n^{\epsilon}$ for some $\epsilon \in (0,1)$, then $k(n) \equiv \epsilon^{-1}$. If instead $(n\rho_{n}) = (\log n)^{c_{2}}$ for some $c_{2}\ge 1$, then $k(n) = \lceil \log n/(c_{2} \log\log n) \rceil$. We remark that all regimes in which $\rho_{n} \rightarrow c_{\rho} > 0$ functionally correspond to the regime where $\rho_{n} \equiv 1$ by appropriate rescaling.

\subsection{First-order approximation (deviations)}
\label{sec:FirstOrder}

Under Assumptions~\ref{assump:evalSignal}~and~\ref{assump:specConcentration}, spectral norm analysis via the Davis-Kahan $\sin \Theta$ theorem \citep[Section~7.3]{bhatia1997matrix} yields that for large $n$ there exists $\bW \equiv \bW_{n} \in\defO_{r}$ such that
\begin{equation}
	\label{eq:spectralNaiveBound}
	\|\hat{\bU}-\bU\bW\| = \probO\left\{(n\rho_{n})^{-1/2}\right\}.
\end{equation}
Equation~(\ref{eq:spectralNaiveBound}) provides a coarse benchmark bound for the quantity $\|\hat{\bU}-\bU\bW\|_{\ttinf}$ (since $\|\cdot\|_{\ttinf} \le \|\cdot\|$), a quantity which is shown below to at times be much smaller.

\begin{theorem}
	\label{thrm:firstOrder}
	Suppose that Assumptions~\ref{assump:sparsity}--\ref{assump:entryConcentration} hold and that $n \rho_{n} = \omega\{(\log n)^{2\xi}\}$ with $r^{1/2} \le (\log n)^{\xi}$. Then there exists $\bW \equiv \bW_{n} \in \defO_{r}$ such that
	\begin{align}
		\label{eq:firstOrder}
		\|\hat{\bU}-\bU\bW\|_{\ttinf}
		&= \probO\left[(n\rho_{n})^{-1/2}
		\times \emph{min}\left\{
		r^{1/2}(\log n)^{\xi}\|\bU\|_{\ttinf}, 1
		\right\}\right].
	\end{align}
\end{theorem}
The bound obtained by two-to-infinity norm methods in Eq.~(\ref{eq:firstOrder}) is demonstrably superior to the bound implied by Eq.~(\ref{eq:spectralNaiveBound}) when $r^{1/2}(\log n)^{\xi}\|\bU\|_{\ttinf} \rightarrow 0$ as $n \rightarrow \infty$, namely when $\|\bU\|_{\ttinf} \rightarrow 0$ sufficiently quickly. Such behavior arises both in theory and in applications, including under the guise of eigenvector delocalization \citep{rudelson2015delocalization,erdHos2013spectral} and of subspace basis coherence \citep{candes2009exact}.

The proof of Theorem~\ref{thrm:firstOrder} first proceeds by way of refined deterministic matrix decompositions and then subsequently leverages the aforementioned probabilistic concentration assumptions. Our proof framework further permits second-order analysis, culminating in Theorem~\ref{thrm:secondOrder} in Section~\ref{sec:SecondOrder}. In the process of proving Theorem~\ref{thrm:secondOrder} we also prove Theorem~\ref{thrm:firstOrderExtension}, an extension and refinement of Theorem~\ref{thrm:firstOrder}. Proof details are provided in the Supplementary Material.

\begin{theorem}
	\label{thrm:firstOrderExtension}
	Suppose that Assumptions~\ref{assump:sparsity}--\ref{assump:entryConcentration} hold and that Eq.~(\ref{eq:entryConcentration}) holds for $k$ up to $k(n)+1$. Suppose $n \rho_{n} = \omega\{(\log n)^{2\xi}\}$ and $r^{1/2} \le (\log n)^{\xi}$. Then there exists $\bW \equiv \bW_{n} \in \defO_{r}$ such that
	\begin{equation}
		\label{eq:firstOrderExtension}
		\hat{\bU}-\bU\bW = \bE\bU\bLam^{-1}\bW + \bR
	\end{equation}
	for some matrix $\bR\in\R^{n \times r}$ satisfying
	\begin{equation*}
		\|\bR\|_{\ttinf}
		= \probO\left[ (n\rho_{n})^{-1} \times r \times \ \emph{max}\left\{(\log n)^{2\xi},\|\bU^{\top}\bE\bU\|+1\right\} \times \|\bU\|_{\ttinf}\right].
	\end{equation*}
	Moreover,
	\begin{equation*}
		\|\bE\bU\bLam^{-1}\bW\|_{\ttinf} = \probO\left\{(n\rho_{n})^{-1/2} \times r^{1/2}(\log n)^{\xi}\|\bU\|_{\ttinf}\right\}.
	\end{equation*}
\end{theorem}

Theorem~\ref{thrm:firstOrderExtension} provides a collective eigenvector (i.e.~subspace) characterization of the relationship between the leading eigenvectors of $\bM$ and $\hat{\bM}$ via the perturbation $\bE$, summarized as
\begin{equation*}
	\hat{\bU}
	\approx \hat{\bM}\bU\bLam^{-1}\bW
	= \bU\bW + \bE\bU\bLam^{-1}\bW.
\end{equation*}
The unperturbed eigenvectors satisfy $\bU\bW \equiv \bM\bU\bLam^{-1}\bW$, leading to the striking observation that the eigenvector perturbation characterization is approximately linear in the perturbation $\bE$.

\begin{remark}
	It always holds that $\|\bU^{\top}\bE\bU\| \le \|\bE\|$, where ``$\le$'' can be replaced by ``$\ll$'' upon invoking Hoeffding-type concentration or more generally \emph{$(C,c,\gamma)$-concentration} \citep{o2018random} for suitable choices of $\bE$. Moreover, $\|\bR\|_{\ttinf} \ll \|\bE\bU\bLam^{-1}\bW\|_{\ttinf}$ holds with high probability in Theorem~\ref{thrm:firstOrderExtension} for numerous regimes in which $n\rho_{n} \rightarrow \infty$ and $\|\bU\|_{\ttinf} \rightarrow 0$.
\end{remark}

\begin{remark}
	\label{rem:moreAssumption2}
	Strictly speaking, Eq.~(\ref{eq:spectralNaiveBound}) holds even when the leading eigenvalues of $\bM$ are not of the same order of magnitude, for the bound is fundamentally given by $C\|\bE\|(|\bLam_{rr}|-\|\bLam_{\perp}\|)^{-1}$. Similarly, the first-order bounds in this paper still hold for $\bLam_{\perp} \neq 0$ provided $\|\bLam_{\perp}\|$ is sufficiently small, in which case na\"{i}ve analysis introduces additional terms of the form $\|\bLam_{\perp}\|\|\bLam^{-1}\|\|\sin\Theta(\hat{\bU},\bU)\|$. In practice the exact spike dimension may be unknown, though it can often be consistently estimated via the ``elbow in the scree plot'' approach \citep{zhu2006automatic} provided $\|\bE\|$ is sufficiently small relative to the leading nonzero eigenvalues of $\bM$.
\end{remark}

\subsection{Second-order limit theory (fluctuations)}
\label{sec:SecondOrder}

This section specifies additional structure on $\bM$ and $\bE$ for the purpose of establishing second-order limit theory. Here, $\bM$ is assumed to have strictly positive leading eigenvalues, reminiscent of a spike covariance or kernel population matrix setting. It is possible though more involved to obtain similar second-order results when $\bM$ is allowed to have both strictly positive and strictly negative leading eigenvalues of the same order. Specifically, such modifications would give rise to considerations involving structured orthogonal matrices and the indefinite orthogonal group.

\begin{assumption}
	\label{assump:errorCLT}
	Suppose that $\bM$ can be written as $\bM \equiv \rho_{n}\bX\bX^{\top} \equiv \bU\bLam\bU^{\top}$ with $\bX=[X_{1}|\dots|X_{n}]^{\top}\in\R^{n \times r}$ and $(n^{-1}\bX^{\top}\bX) \rightarrow \bXi \in \R^{r \times r}$ as $n \rightarrow\infty$ for some symmetric invertible matrix $\bXi$. Also suppose that for a fixed index $i$, the scaled $i$-th row of $\bE\bX$, written as $(n\rho_{n})^{-1/2}(\bE\bX)_{i} = (n\rho_{n})^{-1/2}(\sum_{j=1}^{n}\bE_{ij}X_{j})$, converges in distribution to a centered multivariate normal random vector $Y_{i}\in\R^{r}$ with second moment matrix $\bGamma_{i}\in\R^{r \times r}$.
\end{assumption}
\begin{theorem}
	\label{thrm:secondOrder}
	Suppose that Assumptions~\ref{assump:sparsity}--\ref{assump:errorCLT} hold and that Eq.~(\ref{eq:entryConcentration}) holds for $k$ up to $k(n)+1$. Suppose in addition that $n \rho_{n} = \omega\{(\log n)^{2\xi}\}$, $r^{1/2} \le (\log n)^{\xi}$, and
	\begin{equation}
		\label{eq:secondOrderDecayRequirement}
		\rho_{n}^{-1/2} \times r \times \emph{max}\left\{(\log n)^{2\xi},\|\bU^{\top}\bE\bU\|+1\right\} \times \|\bU\|_{\ttinf} \rightarrow 0
	\end{equation}
	in probability as $n \rightarrow\infty$. Let $\hat{U}_{i}$ and $U_{i}$ be column vectors denoting the $i$-th rows of $\hat{\bU}$ and $\bU$, respectively. Then there exist sequences of orthogonal matrices $(\bW)$ and $(\bW_{\bX})$ depending on $n$ such that the random vector $n\rho_{n}^{1/2}\bW_{\bX}^{\top}(\bW\hat{U}_{i} - U_{i})$ converges in distribution to a centered multivariate normal random vector with covariance matrix $\bSig_{i} = \bXi^{-3/2}\bGamma_{i}\bXi^{-3/2}$, i.e.
	\begin{equation}
		n\rho_{n}^{1/2}\bW_{\bX}^{\top}\left(\bW\hat{U}_{i} - U_{i}\right)
		\Rightarrow
		\mathcal{N}_{r}(\bZero, \bSig_{i}).
	\end{equation}
\end{theorem}
Equation~(\ref{eq:secondOrderDecayRequirement}) amounts to a mild regularity condition that ensures $n\rho_{n}^{1/2}\|\bR\|_{\ttinf} \rightarrow 0$ in probability for $\bR \equiv \bR_{n} \in \R^{n \times r}$ as in Theorem~\ref{thrm:firstOrderExtension}. This condition holds, for example, when $\|\bU\|_{\ttinf} = O\{(\log n)^{c_{3}}n^{-1/2}\}$, in which case the left-hand side of Eq.~(\ref{eq:secondOrderDecayRequirement}) can often be shown to behave as $\probO\{(\log n)^{c_{4}}(n\rho_{n})^{-1/2}\}$ where $(\log n)^{c_{4}}(n\rho_{n})^{-1/2} \rightarrow 0$ as $n \rightarrow \infty$. Such bounds on $\|\bU\|_{\ttinf}$ provably arise when the ratio $(\rmMax_{i}\|X_{i}\|)/(\rmMin_{i}\|X_{i}\|)$ is at most polylogarithmic in $n$.

\begin{remark}[Example:~matrix $\bM$ with kernel-type structure]
	\label{rem:distributionF}
	Let $F$ be a probability distribution defined on $\mathcal{X} \subseteq \R^{r}$, and let $X_{1},\dots,X_{n} \sim F$ be independent random vectors with invertible second moment matrix $\bXi \in \R^{r \times r}$. For $\bX = [X_{1}|\dots|X_{n}]^{\top} \in \R^{n \times r}$, let $\bM = \rho_{n}\bX\bX^{\top} \equiv \bU\bLam\bU^{\top}$, so for each $n$ there exists an $r \times r$ orthogonal matrix $\bW_{\bX}$ such that $\rho_{n}^{1/2}\bX=\bU\bLam^{1/2}\bW_{\bX}$. The strong law of large numbers guarantees that $(n^{-1}\bX^{\top}\bX) \rightarrow \bXi$ almost surely as $n \rightarrow \infty$, and so $\bM$ has $r$ eigenvalues of order $\Theta(n\rho_{n})$ asymptotically almost surely. Moreover, $\|\bU\|_{\ttinf} \le Cn^{-1/2}\|\bX\|_{\ttinf}$ asymptotically almost surely for some constant $C>0$, where $\|\bX\|_{\ttinf}$ can be suitably controlled by imposing additional assumptions, such as taking $\mathcal{X}$ to be bounded or imposing moment assumptions on $\|X_{1}\|$. Conditioning on $\bX$ yields a deterministic choice of $\bM$ for the purposes of Assumption~\ref{assump:errorCLT}.
\end{remark}
\begin{remark}[Example:~matrix $\bE$ and multivariate normality]
	To continue the discussion from Remark~\ref{rem:distributionF},
	let all the entries of $\bE$ be centered and independent up to symmetry with common variance $\sigma_{\bE}^{2}>0$. Then, by the classical multivariate central limit theorem, the asymptotic normality condition in Assumption~\ref{assump:errorCLT} holds and 
	$n\rho_{n}^{1/2}\bW_{\bX}^{\top}(\bW\hat{U}_{i} - U_{i})
		\Rightarrow
		\mathcal{N}_{r}(\bZero,\sigma_{\bE}^{2}\bXi^{-2})$
	by Theorem~\ref{thrm:secondOrder}. There are a variety of other regimes in which the multivariate central limit theorem can be invoked for $(n\rho_{n})^{-1/2}(\sum_{j=1}^{n}\bE_{ij}X_{j})$ in order to satisfy the normality condition in Assumption~\ref{assump:errorCLT}, including when the entries of $\bE$ have heterogeneous variances. In practice, we remark that Assumption~\ref{assump:errorCLT} is structurally milder than Assumption~\ref{assump:entryConcentration} with respect to $\bE$.
\end{remark}

\subsection{Simulations}
\label{sec:Simulation}

The $K$-block stochastic block model \citep{holland1983stochastic} is a simple yet ubiquitous random graph model in which vertices are assigned to one of $K$ possible communities (blocks) and where the adjacency of any two vertices is conditionally independent given the two vertices' community memberships. For stochastic block model graphs on $n$ vertices, the binary symmetric adjacency matrix $\bA \in \{0,1\}^{n \times n}$ can be viewed as an additive perturbation of a (low rank) population edge probability matrix $\bP \in [0,1]^{n \times n}$, $\bA = \bP + \bE$, where for $K$-block model graphs the matrix $\bP$ corresponds to an appropriate dilation of the block edge probability matrix $\bB \in [0,1]^{K \times K}$. In the language of this paper, $\hat{\bM} = \bA$ and $\bM = \bP$. It can be verified that versions of the aforementioned assumptions and hypotheses hold for the following examples. Here we set $\rho_{n} \equiv 1$.

\begin{figure}
	\centering
	\includegraphics[width=0.5\textwidth]{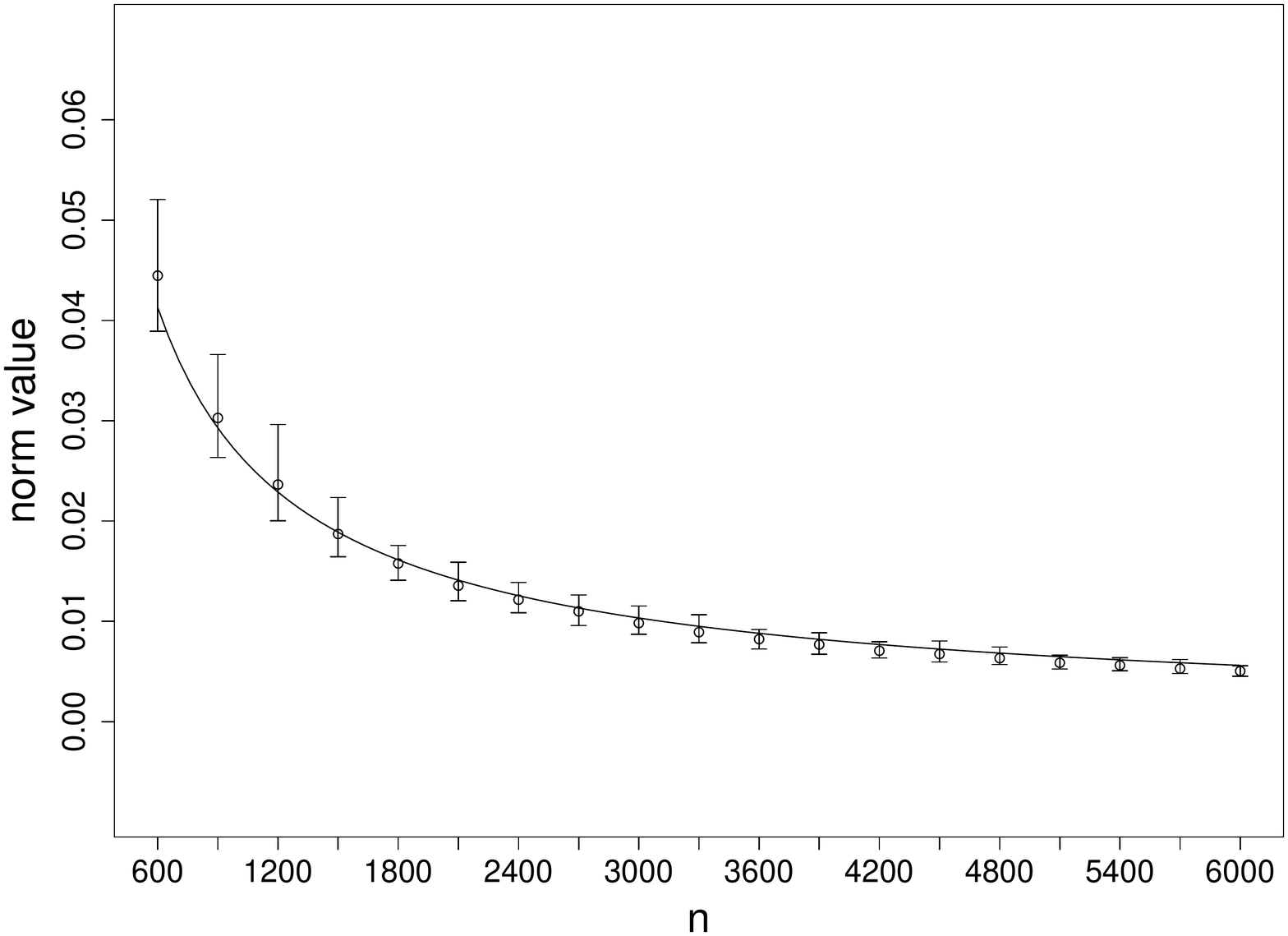}\nolinebreak
	\includegraphics[width=0.5\textwidth]{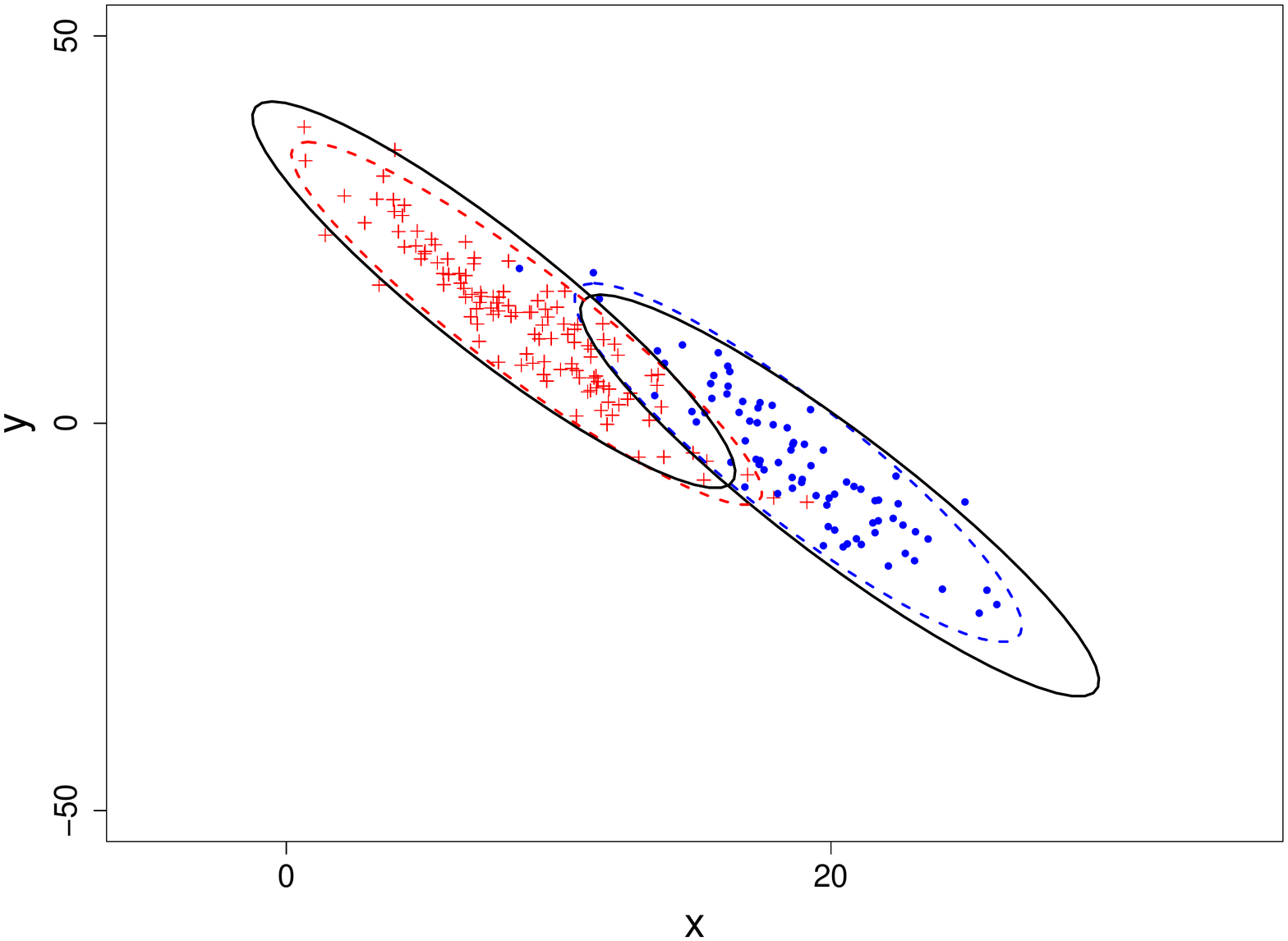}
	\caption{\label{fig:SBMs} (Left plot) First-order simulations for the three-block model with number of vertices $n$ on the $x$-axis and values of $\|\hat{\bU}-\bU\bW\|_{\ttinf}$ on the $y$-axis. Vertical bars depict $95\%$ empirical confidence intervals, and the solid line reflects Theorem~\ref{thrm:firstOrderExtension}. (Right plot) Second-order simulations for the two-block model with $n=200$ where point shape reflects the block membership of the corresponding vertices. Dashed ellipses give the $95\%$ level curves for the empirical distributions. Solid ellipses give the $95\%$ level curves for the theoretical distributions according to Theorem~\ref{thrm:secondOrder}.}
\end{figure}

Consider $n$-vertex graphs arising from the three-block stochastic block model with equal block sizes where the within-block and between-block Bernoulli edge probabilities are given by $\bB_{ii} = 0.5$ for $i=1,2,3$ and $\bB_{ij} = 0.3$ for $i \neq j$, respectively. Here $\rmRank(\bM)=3$, and the second-largest eigenvalue of $\bM$ has multiplicity two. Figure~\ref{fig:SBMs}~(left) plots the empirical mean and $95\%$ empirical confidence interval for $\|\hat{\bU}-\bU\bW\|_{\ttinf}$ computed from $100$ independent simulated adjacency matrices for each value of $n$. Figure~\ref{fig:SBMs}~(left) also plots the function $\phi(n)=\{\lambda_{3}^{-1/2}(\bM)\}(\log n)n^{-1/2}$ which for large $n$ captures the behavior of the leading order term in Theorem~\ref{thrm:firstOrderExtension}. This illustration does not pursue optimal constants or logarithmic factors. Here $\lambda_{3}(\bM)=\Theta(n\rho_{n})=\Theta\{(n\rho_{n})^{1/2}\lambda\}$ with respect to $\lambda$ at the end of Section~\ref{sec:Intro}.

Figure~\ref{fig:SBMs} (right) shows a scatter plot of the (uncentered, block-conditional) scaled leading eigenvector components for an $n=200$ vertex graph arising from a two-block model with $40\%$ of the vertices belonging to the first block and where the block edge probability matrix $\bB$ has entries $\bB_{11} = 0.5$, $\bB_{12} = \bB_{21} = 0.3$, and $\bB_{22} = 0.3$. This small-$n$ example is complemented by additional simulation results provided in the Supplementary Material. We remark that the normalized random (row) vectors are jointly dependent but with decaying pairwise correlations; rows within any fixed finite collection are provably asymptotically independent as $n\rightarrow\infty$.

\section*{Acknowledgment}
We are grateful to the editor, associate editor, and reviewers for their consideration of our paper and for their suggestions. This research was supported by the D3M program of the Defense Advanced Research Projects Agency (DARPA) and by the Acheson J. Duncan Fund for the Advancement of Research in Statistics at Johns Hopkins University.


\section*{Supplementary material}
\label{sec:SM}
	This supplementary section contains a joint proof of the theoretical results in the main paper as well as additional simulation examples.

\subsection{Proofs}
\begin{proof}[of Theorems~2,~3,~and~4]
	We begin with several important observations, namely that
	\begin{equation}
		\label{eq:sinThetaSpectralBound}
		\|(\bI-\bU\bU^{\top})\hat{\bU}\|
		= \|\sin\Theta(\hat{\bU},\bU)\|
		= O\left(\|\bE\||\bLam_{rr}|^{-1}\right)
		= \probO\left\{(n\rho_{n})^{-1/2}\right\},
	\end{equation}
	and that there exists $\bW \in \defO_{r}$ depending on $\hat{\bU}$ and $\bU$ such that
	\begin{equation}
		\label{eq:sinSquaredSpectralBound}
		\|\bU^{\top}\hat{\bU}-\bW\|
		\le \|\sin\Theta(\hat{\bU},\bU)\|^{2}
		= \probO\left\{(n\rho_{n})^{-1}\right\}.
	\end{equation}
	In particular, $\bW$ can be taken to be the product of the left and right orthogonal factors in the singular value decomposition of $\bU^{\top}\hat{\bU}$. Additional details may be found, for example, in \cite{cape2018two}.
	
	Importantly, the relation
	$\hat{\bU}\hat{\bLam} = \hat{\bM}\hat{\bU} = (\bM+\bE)\hat{\bU}$ yields the matrix equation $\hat{\bU}\hat{\bLam}-\bE\hat{\bU}=\bM\hat{\bU}$. The spectra of $\hat{\bLam}$ and $\bE$ are disjoint from one another with high probability as a consequence of Assumptions~2~and~3, so it follows that $\hat{\bU}$ can be written as the matrix series \citep[Section~7.2]{bhatia1997matrix}
	\begin{align}
		\label{eq:vonNeumannExpansion}
		\hat{\bU}
		&= \sum_{k=0}^{\infty}\bE^{k}\bM\hat{\bU}\hat{\bLam}^{-(k+1)}
		= 
		\sum_{k=0}^{\infty}\bE^{k}\bU\bLam\bU^{\top}\hat{\bU}\hat{\bLam}^{-(k+1)},
	\end{align}
	where the second equality holds since $\rmRank(\bM)=r$.
	
	For any choice of $\bW \in \defO_{r}$, the matrix $\hat{\bU}-\bU\bW$ can be decomposed as
	\begin{align*}
		\hat{\bU}-\bU\bW
		&= \bE\hat{\bU}\hat{\bLam}^{-1}
		+ \bU\bLam(\bU^{\top}\hat{\bU}\hat{\bLam}^{-1} - \bLam^{-1}\bU^{\top}\hat{\bU})
		+ \bU(\bU^{\top}\hat{\bU} -\bW)\\
		&= \bE\hat{\bU}\hat{\bLam}^{-1} + \bR^{(1)} + \bR_{\bW}^{(2)}.
	\end{align*}
	For $\bR_{\bW}^{(2)} = \bU(\bU^{\top}\hat{\bU} -\bW)$, it follows that for $\bW$ satisfying Eq.~(\ref{eq:sinSquaredSpectralBound}), then
	\begin{equation*}
		\label{eq:boundR2W}
		\|\bR_{\bW}^{(2)}\|_{\ttinf}
		\le \|\bU^{\top}\hat{\bU}-\bW\|\|\bU\|_{\ttinf}
		= \probO\left\{(n\rho_{n})^{-1}\|\bU\|_{\ttinf}\right\}.
	\end{equation*}
	For $\bR^{(1)} = \bU\bLam\bR^{(3)}$ where $\bR^{(3)} = (\bU^{\top}\hat{\bU}\hat{\bLam}^{-1} - \bLam^{-1}\bU^{\top}\hat{\bU}) \in \R^{r \times r}$, the entries of $\bR^{(3)}$ satisfy
	\begin{equation*}
		\bR_{ij}^{(3)}
		= \langle u_{i}, \hat{u}_{j} \rangle\left\{(\hat{\bLam}_{jj})^{-1} - (\bLam_{ii})^{-1}\right\}
		= \langle u_{i}, \hat{u}_{j} \rangle (\bLam_{ii} - \hat{\bLam}_{jj})
		(\bLam_{ii})^{-1}(\hat{\bLam}_{jj})^{-1}.
	\end{equation*}
	Define the matrix $\bH_{1} \in \R^{r \times r}$ entrywise according to $(\bH_{1})_{ij}=(\bLam_{ii})^{-1}(\hat{\bLam}_{jj})^{-1}$. Then, with $\circ$ denoting the Hadamard matrix product,
	\begin{equation*}
		\bR^{(3)} = - \bH_{1} \circ (\bU^{\top}\hat{\bU}\hat{\bLam} - \bLam\bU^{\top}\hat{\bU}).
	\end{equation*}
	The rightmost matrix factor can be expanded as
	\begin{equation*}
		(\bU^{\top}\hat{\bU}\hat{\bLam} - \bLam\bU^{\top}\hat{\bU})
		= \bU^{\top}\bE\hat{\bU}
		= \bU^{\top}\bE\bU\bU^{\top}\hat{\bU}
		+ \bU^{\top}\bE(\bI-\bU\bU^{\top})\hat{\bU},
	\end{equation*}
	and is therefore bounded in spectral norm using Eq.~(\ref{eq:sinThetaSpectralBound}) in the manner
	\begin{align*}
		\|\bU^{\top}\hat{\bU}\hat{\bLam} - \bLam\bU^{\top}\hat{\bU}\|
		&\le \|\bU^{\top}\bE\bU\|+\probO(1).
	\end{align*}
	Combining the above observations together with properties of matrix norms yields the following two-to-infinity norm bound on $\bR^{(1)}$.
	\begin{align*}
		\|\bR^{(1)}\|_{\ttinf}
		= \|\bU\bLam\bR^{(3)}\|_{\ttinf}
		&\le
		r\|\bU\|_{\ttinf}\|\bLam\|\|\bH_{1}\|_{\rmMax}\|\bU^{\top}\hat{\bU}\hat{\bLam} - \bLam\bU^{\top}\hat{\bU}\|\\
		&=\probO\left\{r(n\rho_{n})^{-1}(\|\bU^{\top}\bE\bU\|+1)\|\bU\|_{\ttinf}\right\}	
	\end{align*}
	
	Assumptions~2~and~3 with an application of Weyl's inequality \citep[Corollary~3.2.6]{bhatia1997matrix} guarantee that there exist constants $C_{1}, C_{2} > 0$ such that $\|\bE\| \le C_{1}(n\rho_{n})^{1/2}$ and $\|\hat{\bLam}^{-1}\| \le C_{2}(n\rho_{n})^{-1}$ with high probability for $n$ sufficiently large. Therefore, by applying the earlier matrix series expansion,
	\begin{align*}
		\|\bE\hat{\bU}\hat{\bLam}^{-1}\|_{\ttinf}
		&= \left\|\sum_{k=1}^{\infty}\bE^{k}\bU\bLam\bU^{\top}\hat{\bU}\hat{\bLam}^{-(k+1)}\right\|_{\ttinf}\\
		&\le \sum_{k=1}^{k(n)}\|\bE^{k}\bU\|_{\ttinf}\|\bLam\|\|\hat{\bLam}^{-1}\|^{k+1}
		+ \sum_{k=k(n)+1}^{\infty}\|\bE\|^{k}\|\bLam\|\|\hat{\bLam}^{-1}\|^{k+1}\\
		&= \probO\left\{r^{1/2}(n\rho_{n})^{-1/2}(\log n)^{\xi}\|\bU\|_{\ttinf} + (n\rho_{n})^{-1/2}\|\bU\|_{\ttinf}\right\},
	\end{align*}
	where we have used the fact that $n\rho_{n} =\omega\{(\log n)^{2\xi}\}$,
	$(n\rho_{n})^{-k(n)/2} \le n^{-1/2} \le \|\bU\|_{\ttinf}$ for $n$ sufficiently large,
	and that by Assumption~4, for each $k \le k(n)$, with high probability
	\begin{equation*}
		\label{eq:highKttinfBound}
		\|\bE^{k}\bU\|_{\ttinf}
		\le r^{1/2}\underset{i\in[n], j\in[r]}{\rmMax}|\langle\bE^{k}u_{j},e_{i}\rangle|
		\le r^{1/2} (C_{\bE}n\rho_{n})^{k/2}(\log n)^{k\xi}\|\bU\|_{\ttinf}.
	\end{equation*}
	Since $\|\bU^{\top}\bE\bU\| \le \|\bE\|$ and
	$ r^{1/2} \le (\log n)^{\xi}$ with $n \rho_{n} = \omega\{(\log n)^{2\xi}\}$, then
	\begin{align*}
		\|\hat{\bU}-\bU\bW\|_{\ttinf}
		&\le \|\bE\hat{\bU}\hat{\bLam}^{-1}\|_{\ttinf} + \|\bR^{(1)}\|_{\ttinf} + \|\bR_{\bW}^{(2)}\|_{\ttinf}\\
		&= \probO\left\{r^{1/2}(n\rho_{n})^{-1/2}(\log n)^{\xi}\|\bU\|_{\ttinf}\right\}.
	\end{align*}
	This completes the proof of Theorem~2.
	
	Next, we further decompose the matrix $\bE\hat{\bU}\hat{\bLam}^{-1}$ by extending the above proof techniques in order to obtain second-order fluctuations. Using the matrix series form in Eq.~(\ref{eq:vonNeumannExpansion}) yields
	\begin{align*}
		\bE\hat{\bU}\hat{\bLam}^{-1}
		&= \bE\bU\bLam\bU^{\top}\hat{\bU}\hat{\bLam}^{-2} + \sum_{k=2}^{\infty}\bE^{k}\bU\bLam\bU^{\top}\hat{\bU}\hat{\bLam}^{-(k+1)}\\
		&= \bE\bU\bLam^{-1}\bW
		+ \bE\bU\bLam(\bU^{\top}\hat{\bU}\hat{\bLam}^{-2}-\bLam^{-2}\bU^{\top}\hat{\bU})
		+ \bE\bU\bLam^{-1}(\bU^{\top}\hat{\bU}-\bW) \\
		&\indent + \sum_{k=2}^{\infty}\bE^{k}\bU\bLam\bU^{\top}\hat{\bU}\hat{\bLam}^{-(k+1)}\\
		&= \bE\bU\bLam^{-1}\bW
		+ \bR_{2}^{(1)} + \bR_{2,\bW}^{(2)} + \bR_{2}^{(\infty)}.
	\end{align*}
	The final term satisfies the bound
	\begin{equation*}
		\|\bR_{2}^{(\infty)}\|_{\ttinf}
		= \probO\left\{r^{1/2}(n\rho_{n})^{-1}(\log n)^{2\xi}\|\bU\|_{\ttinf}\right\},
	\end{equation*}
	which follows from Assumption~4 holding up to $k(n) + 1$, namely
	\begin{align*}
		\|\bR_{2}^{(\infty)}\|_{\ttinf}
		&\le \sum_{k=2}^{k(n)+1}\|\bE^{k}\bU\|_{\ttinf}\|\bLam\|\|\hat{\bLam}^{-1}\|^{k+1} + \sum_{k=k(n)+2}^{\infty}\|\bE\|^{k}\|\bLam\|\|\hat{\bLam}^{-1}\|^{k+1}\\
		&= \probO\left\{r^{1/2}(n\rho_{n})^{-1}(\log n)^{2\xi}\|\bU\|_{\ttinf}+(n \rho_{n})^{-1}\|\bU\|_{\ttinf}\right\}.
	\end{align*}
	On the other hand, modifying the previous analysis used to bound $\bR_{\bW}^{(2)}$ yields
	\begin{equation*}
		\|\bR_{2,\bW}^{(2)}\|_{\ttinf}
		\le \|\bE\bU\|_{\ttinf}\|\bLam^{-1}\|\|\bU^{\top}\hat{\bU}-\bW\|
		= \probO\left\{r^{1/2}(n\rho_{n})^{-3/2}(\log n)^{\xi}\|\bU\|_{\ttinf}\right\}.
	\end{equation*}
	We now bound $\bR_{2}^{(1)}=\bE\bU\bLam(\bU^{\top}\hat{\bU}\hat{\bLam}^{-2}-\bLam^{-2}\bU^{\top}\hat{\bU})$ by extending the previous argument used to bound $\bR^{(1)}$. For $\bR_{2}^{(1)} = \bE\bU\bLam\bR_{2}^{(3)}$ where $\bR_{2}^{(3)} = (\bU^{\top}\hat{\bU}\hat{\bLam}^{-2} - \bLam^{-2}\bU^{\top}\hat{\bU}) \in \R^{r \times r}$, the entries of $\bR_{2}^{(3)}$ satisfy
	\begin{equation*}
		\bR_{ij}^{(3)}
		= \langle u_{i}, \hat{u}_{j} \rangle\left\{(\hat{\bLam}_{jj})^{-2} - (\bLam_{ii})^{-2}\right\}
		= \langle u_{i}, \hat{u}_{j} \rangle (\bLam_{ii}^{2} - \hat{\bLam}_{jj}^{2})
		(\bLam_{ii})^{-2}(\hat{\bLam}_{jj})^{-2}.
	\end{equation*}
	Define the matrix $\bH_{2} \in \R^{r \times r}$ entrywise according to $(\bH_{2})_{ij}=(\bLam_{ii})^{-2}(\hat{\bLam}_{jj})^{-2}$. Then, with $\circ$ denoting the Hadamard matrix product,
	\begin{equation*}
		\bR_{2}^{(3)} = - \bH_{2} \circ (\bU^{\top}\hat{\bU}\hat{\bLam}^{2} - \bLam^{2}\bU^{\top}\hat{\bU}).
	\end{equation*}
	The rightmost matrix factor can be written as
	\begin{equation*}
		(\bU^{\top}\hat{\bU}\hat{\bLam}^{2} - \bLam^{2}\bU^{\top}\hat{\bU})
		= \bU^{\top}(\hat{\bM})^{2}\hat{\bU} - \bU^{\top}\bM^{2}\hat{\bU}
		= \bU^{\top}(\bM\bE+\bE\bM)\hat{\bU},
	\end{equation*}
	and has spectral norm on the order of  $\probO\{(n\rho_{n})^{3/2}\}$. Hence,
	\begin{align*}
		\|\bR_{2}^{(1)}\|_{\ttinf}
		= \|\bE\bU\bLam\bR_{2}^{(3)}\|_{\ttinf}
		&\le
		r\|\bE\bU\|_{\ttinf}\|\bLam\|\|\bH_{2}\|_{\rmMax}\|\bU^{\top}\hat{\bU}\hat{\bLam}^{2} - \bLam^{2}\bU^{\top}\hat{\bU}\|\\
		&= \probO\left\{r^{3/2}(n\rho_{n})^{-1}(\log n)^{\xi}\|\bU\|_{\ttinf}\right\}.
	\end{align*}	
	For $\bR = \bR^{(1)} + \bR_{\bW}^{(2)} + \bR_{2}^{(1)} +\bR_{2,\bW}^{(2)} + \bR_{2}^{(\infty)}$, we have therefore shown that
	\begin{equation}
		\label{eq:decompositionOverview}
		\hat{\bU}-\bU\bW = \bE\bU\bLam^{-1}\bW + \bR,
	\end{equation}
	where since $r^{1/2} \le (\log n)^{\xi}$, the residual matrix $\bR$ satisfies
	\begin{equation*}
		\|\bR\|_{\ttinf}
		= \probO\left[(n\rho_{n})^{-1} \times r \times \rmMax\left\{(\log n)^{2\xi},\|\bU^{\top}\bE\bU\|+1\right\} \times \|\bU\|_{\ttinf}\right].
	\end{equation*}
	The leading term agrees with the order of the bound in Theorem~2, namely
	\begin{equation*}
		\|\bE\bU\bLam^{-1}\bW\|_{\ttinf} = \probO\left\{(n\rho_{n})^{-1/2} \times r^{1/2}(\log n)^{\xi}\|\bU\|_{\ttinf}\right\}.
	\end{equation*}
	This establishes Theorem~3 en route to proving Theorem~4, which we now proceed to finish.
	
	Since $\bM = \rho_{n}\bX\bX^{\top} \equiv \bU\bLam\bU^{\top}$, there exists an orthogonal matrix $\bW_{\bX}$ (depending on $n$) such that $\rho_{n}^{1/2}\bX=\bU\bLam^{1/2}\bW_{\bX}$, hence $\rho_{n}\bX^{\top}\bX=\bW_{\bX}^{\top}\bLam\bW_{\bX}$. Following some algebraic manipulations, the matrix $\bE\bU\bLam^{-1}\bW$ can therefore be written as
	\begin{align*}
		\bE\bU\bLam^{-1}\bW
		&= \rho_{n}^{-1}\bE\bX(\bX^{\top}\bX)^{-3/2}(\bW_{\bX}^{\top}\bW).
	\end{align*}
	Plugging this observation into Eq.~(\ref{eq:decompositionOverview}) and subsequent matrix multiplication together yield the relation
	\begin{align*}
		\left(\hat{\bU}\bW^{\top}\bW_{\bX}-\bU\bW_{\bX}\right)
		= \rho_{n}^{-1}\bE\bX(\bX^{\top}\bX)^{-3/2}
		+ \bR\bW^{\top}\bW_{\bX}.
	\end{align*}
	For fixed $i$, let $\hat{U}_{i}$, $U_{i},$ and $R_{i}$ be column vectors denoting the $i$-th rows of $\hat{\bU}$, $\bU$, and $\bR$, respectively. Equation~(7) in the main paper implies that $n\rho_{n}^{1/2}\|R_{i}\| \rightarrow 0$ in probability. In addition,  $(n^{-1}\bX^{\top}\bX)^{-3/2} \rightarrow \bXi^{-3/2}$ by Assumption~5 together with the continuous mapping theorem. The scaled $i$-th row of $\bE\bX$ converges in distribution to $Y_{i} \sim \mathcal{N}_{r}(\bZero, \bGamma_{i})$ by Assumption~5, so combining the above observations together with Slutsky's theorem yields that there exist sequences of orthogonal matrices $(\bW)$ and $(\bW_{\bX})$ such that
	\begin{align*}
		n\rho_{n}^{1/2}\bW_{\bX}^{\top}\left(\bW\hat{U}_{i} - U_{i}\right)
		&\overset{}{=} \left(n^{-1}\bX^{\top}\bX\right)^{-3/2}\left\{(n\rho_{n})^{-1/2}(\bE\bX)_{i}\right\} + n\rho_{n}^{1/2}\bW_{\bX}^{\top}\bW R_{i}\\
		&\Rightarrow
		\bXi^{-3/2}Y_{i} + 0.
	\end{align*}
	In particular, we have the row-wise convergence in distribution
	\begin{equation*}
		n\rho_{n}^{1/2}\bW_{\bX}^{\top}\left(\bW\hat{U}_{i} - U_{i}\right)
		\Rightarrow
		\mathcal{N}_{r}(\bZero,\bSig_{i})
	\end{equation*}
	where $\bSig_{i} = \bXi^{-3/2}\bGamma_{i}\bXi^{-3/2}$. This completes the proof of Theorem~4.
\end{proof}

\subsection{Two-block stochastic block model (continued)}
Consider $n$-vertex graphs arising from the two-block stochastic block model with $40\%$ of the vertices belonging to the first block and where the block edge probability matrix $\bB$ has entries $\bB_{11} = 0.5$, $\bB_{12} = \bB_{21} = 0.3$, and $\bB_{22} = 0.3$. This model corresponds to Figure~1 (right) in the main paper.
Here, Table~\ref{fig:secondOrderTable} shows block-conditional sample covariance matrix estimates for the centered random vectors $n \rho_{n}^{1/2}\bW_{\bX}^{\top}(\bW\hat{U}_{i} - U_{i})$. Also shown are the corresponding theoretical covariance matrices.

\begin{table}[!]
	\begin{center}
		\caption{\label{fig:secondOrderTable} Empirical and theoretical covariance matrices for the two-block model}
		\begin{tabular}{c c c c} 
			\hline
			$n$ & 1000 & 2000 & $\infty$ \\ [0.5ex] 
			\hline\hline
			\vspace{1em}
			$\hat{\bSig}_{1}$ &
			$\left[\begin{array}{r r}
			14.11 & -36.08 \\
			-36.08 & 110.13
			\end{array}\right]$ &
			$\left[\begin{array}{r r}
			14.94 & -36.85 \\
			-36.85 & 108.55
			\end{array}\right]$ &
			$\left[\begin{array}{r r}
			15.14 & -38.05 \\
			-38.05 & 112.34
			\end{array}\right]$ \\
			$\hat{\bSig}_{2}$ &
			$\left[\begin{array}{r r}
			11.76 & -30.09 \\
			-30.09 & 93.07
			\end{array}\right]$ &
			$\left[\begin{array}{r r}
			12.91 & -33.04 \\
			-33.04 & 101.64
			\end{array}\right]$ &
			$\left[\begin{array}{r r}
			13.12 & -33.93 \\
			-33.93 & 103.94
			\end{array}\right]$
		\end{tabular}
	\end{center}
\end{table}

\subsection{Spike matrix models}
\begin{figure}[!]
	\centering
	\includegraphics[width=0.5\textwidth]{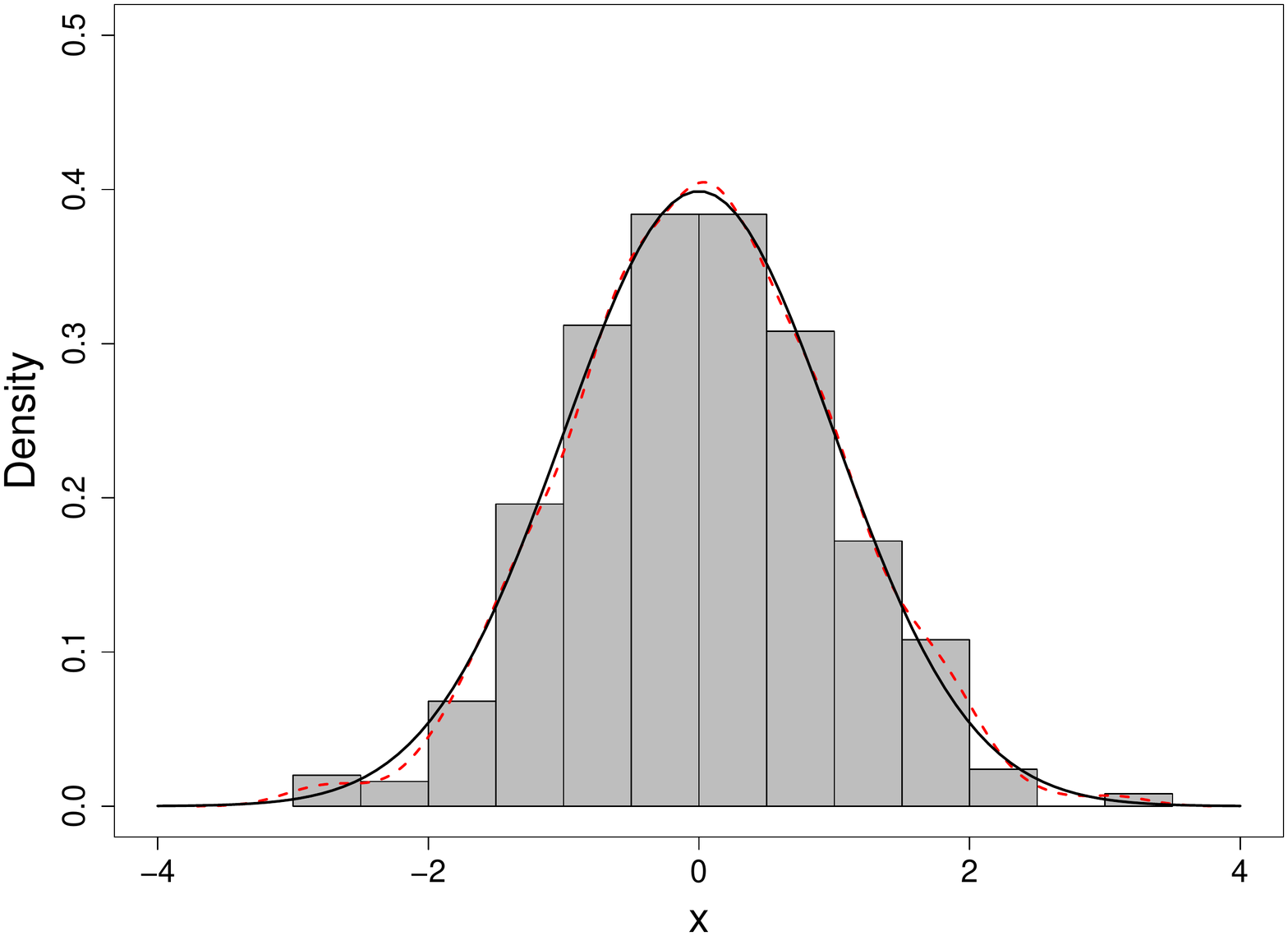}\nolinebreak
	\includegraphics[width=0.5\textwidth]{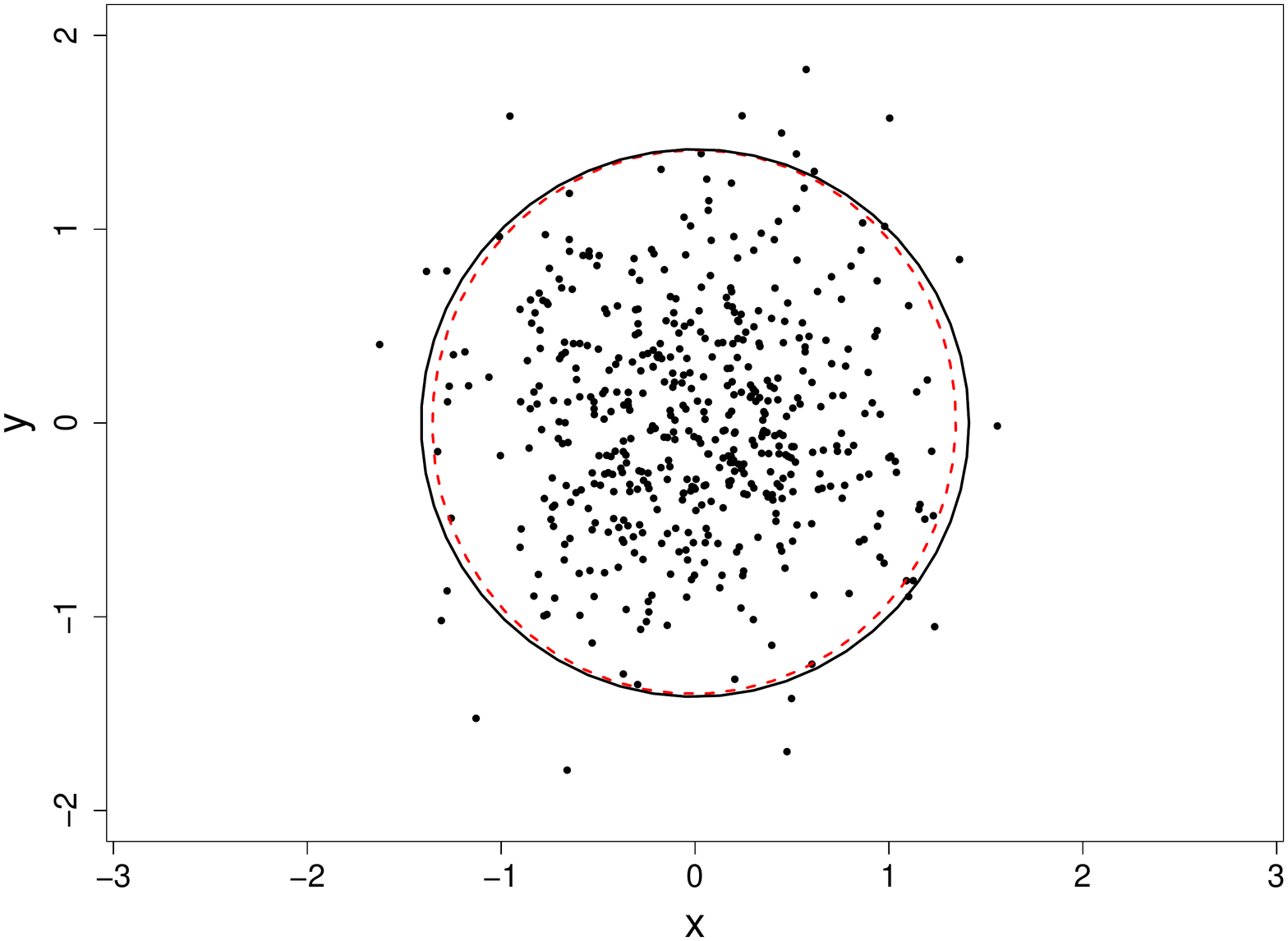}
	\caption{\label{fig:spikePlots} (Left plot) One-dimensional simulation for $n=500$ with empirical (dashed line) and theoretical (solid line) eigenvector fluctuation density. (Right plot) Two-dimensional simulation for $n=500$ where the dashed ellipse gives the $95\%$ level curve for the empirical distribution, and the solid ellipse gives the $95\%$ level curve for the row-wise theoretical distribution.}
\end{figure}

Figure~\ref{fig:spikePlots} provides two additional examples illustrating Theorem~4 in the main paper for one and two-dimensional spike matrix models, written in the rescaled form $\hat{\bM} = \lambda \bU\bU^{\top} + \bE$ with $\rho_{n} \equiv 1$.
In the left plot, $\lambda=n$, $\bU=n^{-1/2}e\in\R^{n}$, and $\bE_{ij} \sim \textnormal{Laplace}(0,2^{-1/2})$ independently for $i \le j$ with $\bE_{ij} = \bE_{ji}$. Here $\bXi$ is the one-dimensional identity matrix, i.e.~$\bXi = \bI_{1}$, and $(n\rho_{n})^{-1/2}(\bE\bX)_{i} \Rightarrow \mathcal{N}_{1}(0,1)$ by the central limit theorem, so for each fixed row $i$ Theorem~4 yields convergence in distribution to $\mathcal{N}_{1}(0,1)$.
In the right plot, $\lambda = n$ and $\bU_{ij}=n^{-1/2}$ for $1\le i \le n, j=1$, $1\le i \le n/2, j=2$ with $\bU_{ij}=-n^{-1/2}$ otherwise. In addition, $\bE_{ij} \sim \textnormal{Uniform}[-1,1]$ independently for $i \le j$ with $\bE_{ij} = \bE_{ji}$, so $\textnormal{Var}(\bE_{ij})=1/3$. 
Here $(n\rho_{n})^{-1/2}(EX)_{i}$ converges in distribution to a centered multivariate normal random variable with covariance matrix $\bGamma_{i}=(1/3)\bI_{2} \in \R^{2 \times 2}$ by the multivariate central limit theorem, while the second moment matrix for the rows $X_{i}$ in Assumption~5 is simply $\Xi=\bI_{2}$. Theorem~4 therefore yields 
$n\rho_{n}^{1/2}\bW_{\bX}^{\top}(\bW\hat{U}_{i} - U_{i})
\Rightarrow
\mathcal{N}_{2}(\mu, \Sigma_{i})$,
where $\mu = (0,0)^{\top} \in \R^{2}$ and $\Sigma_{i} = (1/3)I_{2} \in \R^{2 \times 2}$.
Plots depict all vectors computed from a single simulated adjacency matrix.

\bibliographystyle{biometrika}
\bibliography{unscaled}

\begin{thebibliography}{}

\bibitem[\protect\citeauthoryear{Abbe, Fan, Wang, and Zhong}{Abbe
  et~al.}{2017}]{abbe2017entrywise}
Abbe, E., J.~Fan, K.~Wang, and Y.~Zhong (2017).
\newblock Entrywise eigenvector analysis of random matrices with low expected
  rank.
\newblock {\em preprint arXiv:1709.09565\/}.

\bibitem[\protect\citeauthoryear{Bai and Silverstein}{Bai and
  Silverstein}{2010}]{bai2010spectral}
Bai, Z. and J.~W. Silverstein (2010).
\newblock {\em Spectral analysis of large dimensional random matrices},
  Volume~20.
\newblock Springer.

\bibitem[\protect\citeauthoryear{Benaych-Georges and
  Nadakuditi}{Benaych-Georges and Nadakuditi}{2011}]{benaych2011eigenvalues}
Benaych-Georges, F. and R.~R. Nadakuditi (2011).
\newblock The eigenvalues and eigenvectors of finite, low rank perturbations of
  large random matrices.
\newblock {\em Advances in Mathematics\/}~{\em 227\/}(1), 494--521.

\bibitem[\protect\citeauthoryear{Bhatia}{Bhatia}{1997}]{bhatia1997matrix}
Bhatia, R. (1997).
\newblock {\em Matrix Analysis}, Volume 169 of Graduate Texts in Mathematics.
\newblock Springer-Verlag, New York.

\bibitem[\protect\citeauthoryear{Cai and Zhang}{Cai and
  Zhang}{2018}]{cai2018rate}
Cai, T.~T. and A.~Zhang (2018).
\newblock Rate-optimal perturbation bounds for singular subspaces with
  applications to high-dimensional statistics.
\newblock {\em The Annals of Statistics\/}~{\em 46\/}(1), 60--89.

\bibitem[\protect\citeauthoryear{Cand{\`e}s and Recht}{Cand{\`e}s and
  Recht}{2009}]{candes2009exact}
Cand{\`e}s, E.~J. and B.~Recht (2009).
\newblock Exact matrix completion via convex optimization.
\newblock {\em Foundations of Computational Mathematics\/}~{\em 9\/}(6), 717.

\bibitem[\protect\citeauthoryear{Cape, Tang, and Priebe}{Cape
  et~al.}{2018}]{cape2018two}
Cape, J., M.~Tang, and C.~E. Priebe (2018).
\newblock The two-to-infinity norm and singular subspace geometry with
  applications to high-dimensional statistics.
\newblock {\em The Annals of Statistics, accepted, preprint
  arXiv:1705.10735\/}.

\bibitem[\protect\citeauthoryear{Eldridge, Belkin, and Wang}{Eldridge
  et~al.}{2018}]{eldridge2018unperturbed}
Eldridge, J., M.~Belkin, and Y.~Wang (2018).
\newblock Unperturbed: spectral analysis beyond {D}avis-{K}ahan.
\newblock In {\em Proceedings of Algorithmic Learning Theory}, Volume~83 of
  {\em Proceedings of Machine Learning Research}, pp.\  321--358. PMLR.

\bibitem[\protect\citeauthoryear{Erd{\H{o}}s, Knowles, Yau, and
  Yin}{Erd{\H{o}}s et~al.}{2013}]{erdHos2013spectral}
Erd{\H{o}}s, L., A.~Knowles, H.-T. Yau, and J.~Yin (2013).
\newblock Spectral statistics of {E}rd{\H{o}}s--{R}{\'e}nyi graphs {I}: {L}ocal
  semicircle law.
\newblock {\em The Annals of Probability\/}~{\em 41\/}(3B), 2279--2375.

\bibitem[\protect\citeauthoryear{Fan, Wang, and Zhong}{Fan
  et~al.}{2018}]{fan2018ell}
Fan, J., W.~Wang, and Y.~Zhong (2018).
\newblock An {$\ell_{\infty}$} eigenvector perturbation bound and its
  application to robust covariance estimation.
\newblock {\em Journal of Machine Learning Research\/}~{\em 18\/}(207), 1--42.

\bibitem[\protect\citeauthoryear{Holland, Laskey, and Leinhardt}{Holland
  et~al.}{1983}]{holland1983stochastic}
Holland, P.~W., K.~B. Laskey, and S.~Leinhardt (1983).
\newblock Stochastic blockmodels: {F}irst steps.
\newblock {\em Social Networks\/}~{\em 5\/}(2), 109--137.

\bibitem[\protect\citeauthoryear{Johnstone}{Johnstone}{2001}]{johnstone2001PCA}
Johnstone, I.~M. (2001).
\newblock On the distribution of the largest eigenvalue in principal components
  analysis.
\newblock {\em The Annals of Statistics\/}~{\em 29\/}(2), 295--327.

\bibitem[\protect\citeauthoryear{Jolliffe}{Jolliffe}{1986}]{jolliffe1986principal}
Jolliffe, I.~T. (1986).
\newblock {\em {P}rincipal {C}omponent {A}nalysis}.
\newblock Springer.

\bibitem[\protect\citeauthoryear{Le, Levina, and Vershynin}{Le
  et~al.}{2017}]{le2017concentration}
Le, C.~M., E.~Levina, and R.~Vershynin (2017).
\newblock Concentration and regularization of random graphs.
\newblock {\em Random Structures \& Algorithms\/}~{\em 51\/}(3), 538--561.

\bibitem[\protect\citeauthoryear{Lei and Rinaldo}{Lei and
  Rinaldo}{2015}]{lei2015consistency}
Lei, J. and A.~Rinaldo (2015).
\newblock Consistency of spectral clustering in stochastic block models.
\newblock {\em The Annals of Statistics\/}~{\em 43\/}(1), 215--237.

\bibitem[\protect\citeauthoryear{Mao, Sarkar, and Chakrabarti}{Mao
  et~al.}{2017}]{mao2017estimating}
Mao, X., P.~Sarkar, and D.~Chakrabarti (2017).
\newblock Estimating mixed memberships with sharp eigenvector deviations.
\newblock {\em preprint arXiv:1709.00407\/}.

\bibitem[\protect\citeauthoryear{Nadler}{Nadler}{2008}]{nadler2008finite}
Nadler, B. (2008).
\newblock Finite sample approximation results for principal component analysis:
  A matrix perturbation approach.
\newblock {\em The Annals of Statistics\/}~{\em 36\/}(6), 2791--2817.

\bibitem[\protect\citeauthoryear{O'Rourke, Vu, and Wang}{O'Rourke
  et~al.}{2018}]{o2018random}
O'Rourke, S., V.~Vu, and K.~Wang (2018).
\newblock Random perturbation of low rank matrices: Improving classical bounds.
\newblock {\em Linear Algebra and its Applications\/}~{\em 540}, 26--59.

\bibitem[\protect\citeauthoryear{Paul}{Paul}{2007}]{paul2007asymptotics}
Paul, D. (2007).
\newblock Asymptotics of sample eigenstructure for a large dimensional spiked
  covariance model.
\newblock {\em Statistica Sinica\/}~{\em 17\/}(4), 1617--1642.

\bibitem[\protect\citeauthoryear{Paul and Aue}{Paul and
  Aue}{2014}]{paul2014random}
Paul, D. and A.~Aue (2014).
\newblock Random matrix theory in statistics: A review.
\newblock {\em Journal of Statistical Planning and Inference\/}~{\em 150},
  1--29.

\bibitem[\protect\citeauthoryear{Rohe, Chatterjee, and Yu}{Rohe
  et~al.}{2011}]{rohe2011spectral}
Rohe, K., S.~Chatterjee, and B.~Yu (2011).
\newblock Spectral clustering and the high-dimensional stochastic blockmodel.
\newblock {\em The Annals of Statistics\/}~{\em 39\/}(4), 1878--1915.

\bibitem[\protect\citeauthoryear{Rudelson and Vershynin}{Rudelson and
  Vershynin}{2015}]{rudelson2015delocalization}
Rudelson, M. and R.~Vershynin (2015).
\newblock Delocalization of eigenvectors of random matrices with independent
  entries.
\newblock {\em Duke Mathematical Journal\/}~{\em 164\/}(13), 2507--2538.

\bibitem[\protect\citeauthoryear{Sarkar and Bickel}{Sarkar and
  Bickel}{2015}]{sarkar2015role}
Sarkar, P. and P.~J. Bickel (2015).
\newblock Role of normalization in spectral clustering for stochastic
  blockmodels.
\newblock {\em The Annals of Statistics\/}~{\em 43\/}(3), 962--990.

\bibitem[\protect\citeauthoryear{Silverstein}{Silverstein}{1984}]{silverstein1984some}
Silverstein, J.~W. (1984).
\newblock Some limit theorems on the eigenvectors of large dimensional sample
  covariance matrices.
\newblock {\em Journal of Multivariate Analysis\/}~{\em 15\/}(3), 295--324.

\bibitem[\protect\citeauthoryear{Silverstein}{Silverstein}{1989}]{silverstein1989eigenvectors}
Silverstein, J.~W. (1989).
\newblock On the eigenvectors of large dimensional sample covariance matrices.
\newblock {\em Journal of Multivariate Analysis\/}~{\em 30\/}(1), 1--16.

\bibitem[\protect\citeauthoryear{Tang, Cape, and Priebe}{Tang
  et~al.}{2017}]{tang2017asymptotically}
Tang, M., J.~Cape, and C.~E. Priebe (2017).
\newblock Asymptotically efficient estimators for stochastic blockmodels: the
  naive {MLE}, the rank-constrained {MLE}, and the spectral.
\newblock {\em preprint arXiv:1710.10936\/}.

\bibitem[\protect\citeauthoryear{Tang and Priebe}{Tang and
  Priebe}{2018}]{tang2018limit}
Tang, M. and C.~E. Priebe (2018).
\newblock Limit theorems for eigenvectors of the normalized {L}aplacian for
  random graphs.
\newblock {\em The Annals of Statistics\/}~{\em 46\/}(5), 2360--2415.

\bibitem[\protect\citeauthoryear{Yu, Wang, and Samworth}{Yu
  et~al.}{2014}]{yu2014useful}
Yu, Y., T.~Wang, and R.~J. Samworth (2014).
\newblock A useful variant of the {D}avis--{K}ahan theorem for statisticians.
\newblock {\em Biometrika\/}~{\em 102\/}(2), 315--323.

\bibitem[\protect\citeauthoryear{Zhu and Ghodsi}{Zhu and
  Ghodsi}{2006}]{zhu2006automatic}
Zhu, M. and A.~Ghodsi (2006).
\newblock Automatic dimensionality selection from the scree plot via the use of
  profile likelihood.
\newblock {\em Computational Statistics \& Data Analysis\/}~{\em 51\/}(2),
  918--930.

\end{thebibliography}


\end{document}